\RequirePackage{ifpdf}
\ifpdf % We are running pdfTeX in pdf mode
\documentclass[pdftex]{sigma}
\else
\documentclass{sigma}
\fi

\def\dA{\hbox{$\mid \! \! \! \! \square$}}
\begin{document}

\allowdisplaybreaks
	
\renewcommand{\PaperNumber}{102}

\FirstPageHeading

\renewcommand{\thefootnote}{$\star$}

\ShortArticleName{Translation to Bundle Operators}

\ArticleName{Translation to Bundle Operators\footnote{This paper is a
contribution to the Proceedings of the 2007 Midwest
Geometry Conference in honor of Thomas~P.\ Branson. The full collection is available at
\href{http://www.emis.de/journals/SIGMA/MGC2007.html}{http://www.emis.de/journals/SIGMA/MGC2007.html}}}

\Author{Thomas P. BRANSON~$^\dag$ and Doojin HONG~$^\ddag$}

\AuthorNameForHeading{T. Branson and D. Hong}

\Address{$^\dag$~Deceased}
\URLaddressD{\url{http://www.math.uiowa.edu/~branson/}}

\Address{$^\ddag$~Department of Mathematics, University of North Dakota, Grand Forks ND 58202, USA}
\EmailD{\href{mailto:doojin.hong@und.edu}{doojin.hong@und.edu}}

\ArticleDates{Received August 31, 2007, in f\/inal form October 24, 2007; Published online October 31, 2007}

\Abstract{We give explicit formulas for conformally invariant
operators with leading term an $m$-th power of Laplacian on the
product of spheres with the natural pseudo-Riemannian product metric for all $m$.}


\Classification{53A30; 53C50}

\section{Introduction}

Conformally invariant operators have been one of the major subjects in mathematics and physics.
Getting explicit formulas of such operators on many manifolds is potentially important.
One use of spectral data, among other things, would be in application to
Polyakov formulas in even dimensions for the quotient of functional determinants
of operators since the precise form of these Polyakov formulas only depends on some
constants that appear in the spectral asymptotics of the operators in question \cite{Branson:95}.

In 1987, Branson \cite{Branson:87} showed explicit formulas of invariant operators on functions and dif\/fe\-ren\-tial
forms over the double cover $S^1\times S^{n-1}$ of the $n$ dimensional compactif\/ied
Minkowski space. And lately, Branson and Hong \cite{BH:06, Hong:04, Hong:00} gave explicit determinant quotient formulas
of operators on spinors and twistors including the Dirac and Rarita Schwinger operators over
$S^1\times S^{n-1}$.
Gover \cite{Gover:06} recently exhibited explicit formulas of invariant operators with leading term a
power of Laplacian on functions over conformally Einstein manifolds.

In this paper, we show explicit formulas of invariant operators with leading term a power of
Laplacian on functions over general product of spheres, $S^p\times S^q$ with the natural
pseudo-Riemannian metric.

\section{Yamabe and Paneitz operators}

Consider $S^p\times S^q$ with the natural signature
$(p, q)$ metric ($p$ minus signs), with $p+q=n$.
We view this as imbedded in the natural way in ${\mathbb R}^{n+2}$,
which carries a signature $(p+1,q+1)$ metric~\cite{Branson:87, Orsted:81}
denoting this manifold with metric by ${\mathbb R}^{p+1,q+1}$.
We consider the radial vector f\/ields
\begin{gather*}
s\partial_s=S=x^a\partial_a  \quad \mbox{in the ambient }{\mathbb R}^{p+1}, \\
r\partial_r=R=x^b\partial_b  \quad \mbox{in the ambient }{\mathbb R}^{q+1}.
\end{gather*}
The d'Alembertian $\dA$ on ${\mathbb R}^{p+1,q+1}$ is
\begin{gather}
\dA=\triangle_{{\mathbb R}^{q+1}}-\triangle_{{\mathbb R}^{p+1}}=-\partial^2_r-\frac{q}{r}\partial_r
+r^{-2}\triangle_{S^q}
+\partial^2_s+\frac{p}{s}\partial_s-s^{-2}\triangle_{S^p} \nonumber\\
\phantom{\dA}{} =r^{-2}\left(-R^2-(q-1)R+\triangle_{S^q}\right)-
s^{-2}\left(-S^2-(p-1)S+\triangle_{S^p}\right) ,\label{dAlem}
\end{gather}
where $\Delta=-g^{ab}\nabla_a\nabla_b$.

It is well known that the following process is conformally
invariant \cite{FG:85, HH:83}:
\begin{itemize}\itemsep=0pt
\item Take a function $f$ on $S^p\times S^q$, and extend
it to a function $F$ having
\begin{gather}\label{Xext}
XF=\left(m- {\frac{n}2}\right)F,
\end{gather}
where
\begin{displaymath}
X:=R+S.
\end{displaymath}
\item Compute $\dA^mF$.
\item Restrict to $S^p\times S^q$.
\end{itemize}
More precisely, if we view $f$ as an $(m-n/2)$-density on the product
of spheres, perform the process above, and view the restricted function
as a $-(m+n/2)$-density, we get a conformally invariant operator
\[
{\cal E}[m-n/2]\to{\cal E}[-m-n/2],
\]
where ${\cal E}[\omega]$ is the bundle of conformal densities of degree $\omega$ \cite{GJMS:92}:
\[
f\in{\cal E}[\omega]\iff \hat{f}=\Omega^{\omega}f \mbox{ under }
\hat{g}=\Omega^2g,\quad \Omega\mbox{ is a positive smooth function}.
\]
In fact, this happens in the more general setting of
the Fef\/ferman--Graham ambient space
for a pseudo-Riemannian conformal manifold $(M,[g])$, provided the dimension
is odd, $2m\le n$, or the Fef\/ferman--Graham obstruction tensor vanishes \cite{FG:85}.
In particular, this happens with no restriction on $(n,m)$ whenever
$[g]$ is a f\/lat conformal structure and this is the case in our situation.
In particular, using only invariance under conformal changes
implemented by dif\/feomorphisms, in our situation
we get an intertwining operator $A$ for two representations
of the conformal group O$(p+1,q+1)$ \cite{Branson:87, Branson:96}:
\[
Au_{m-n/2}=u_{-m-n/2}A.
\]
We begin with a function $f$
having homogeneity $u$ in the radial ($S$) direction
in
${\mathbb R}^{p+1}$, and homogeneity $v$ in the radial ($R$) direction in
${\mathbb R}^{q+1}$.

The $(u,v)$ homogeneity extension is a special case of
the extension scheme (\ref{Xext}) as long as
\begin{gather}\label{correctw}
\omega:=u+v=m-n/2.
\end{gather}
To illustrate our method, we work out the Yamabe operator~($m=1$) and
the Paneitz operator~($m=2$) cases.

Let $Y:=R-S$. On $S^p\times S^q$, $r=s=1$ and we have
\begin{gather*}
\dA f\mid_{S^p\times S^q}=\left\{-R^2+S^2-(q-1)R+(p-1)S+\triangle_{S^q}-\triangle_{S^p}\right\}f\\
\phantom{\dA f\mid_{S^p\times S^q}}{}  =\left\{\left(\begin{array}{c} -XY\\ {\rm or}\\-YX\end{array}\right)
   +X\underbrace{\left\{-\frac{q-1}{2}+\frac{p-1}{2}\right\}}_{-\frac{q-p}{2}}
   +Y\underbrace{\left\{-\frac{q-1}{2}-\frac{p-1}{2}\right\}}_{-\frac{n}{2}+1}
   +\square_{S^p\times S^q}\right\}f\\
\phantom{\dA f\mid_{S^p\times S^q}}{}  = {-\omega Yf-\frac{q-p}{2}\omega f+\left(-\frac{n}{2}+1\right)Yf+
   \square_{S^p\times S^q}}f , \qquad \mbox{since} \ \ Xf=\omega f.
\end{gather*}
Thus, if we choose $\omega=1-\frac{n}{2}$, we get
\[
\dA \mid_{S^p\times S^q}=\square_{S^p\times S^q}-\frac{q-p}{2}\left(1-\frac{n}{2}\right).
\]
The scalar curvature on $S^p\times S^q$ is
\[
q(q-1)-p(p-1)=(q+p)(q-p)-(q-p)=(n-1)(q-p)
\]
and we get the Yamabe operator
\[
\dA \mid_{S^p\times S^q}=\square_{S^p\times S^q}+\frac{n-2}{4(n-1)}\mbox{Scal} .
\]

Now we look at $\dA^2$ in the ambient space.
Since $R$, $S$, $\triangle_{S^p}$, and $\triangle_{S^q}$ all commute,
\begin{gather*}
 \dA^2 = r^{-4}\triangle_{S^q}^2+s^{-4}\triangle_{S^p}^2-2r^{-2}s^{-2}\triangle_{S^q}\triangle_{S^p} \\
\phantom{\dA^2 =}{}  +r^{-2}\triangle_{S^q}\left(r^{-2}\{-R^2-(q-1)R\}+s^{-2}\{S^2+(p-1)S\}\right) \\
\phantom{\dA^2 =}{}  +\left(r^{-2}\{-R^2-(q-1)R\}+s^{-2}\{S^2+(p-1)S\}\right)(r^{-2}\triangle_{S^q}) \\
\phantom{\dA^2 =}{}  -s^{-2}\triangle_{S^p}\left(r^{-2}\{-R^2-(q-1)R\}+s^{-2}\{S^2+(p-1)S\}\right) \\
\phantom{\dA^2 =}{}  -\left(r^{-2}\{-R^2-(q-1)R\}+s^{-2}\{S^2+(p-1)S\}\right)(s^{-2}\triangle_{S^p}) \\
\phantom{\dA^2 =}{} +r^{-2}\{-R^2-(q-1)R\}\left(r^{-2}\{-R^2-(q-1)R\}\right) \\
\phantom{\dA^2 =}{} +s^{-2}\{S^2+(p-1)S\}\left(s^{-2}\{S^2+(p-1)S\}\right) \\
\phantom{\dA^2 =}{}  +2r^{-2}s^{-2}\{-R^2-(q-1)R\}\{S^2+(p-1)S\} .
\end{gather*}
Since $r^{-2}\{-R^2-(q-1)R\}\{r^{-2}(-R^2-(q-1)R)\}$ equals
\begin{gather*}
 r^{-2}R\left(r^{-2}R\{R^2+(q-1)R\}-2r^{-2}\{R^2+(q-1)R\}\right) \\
 \qquad{}+r^{-2}(q-1)\left(r^{-2}R\{R^2+(q-1)R\}-2r^{-2}\{R^2+(q-1)R\}\right) \\
\qquad{} = r^{-4}\left((R^2-4R+4)\{R^2+(q-1)R\}+(q-1)(R-2)\{R^2+(q-1)R\}\right) \\
 \qquad{}= r^{-4}\{R^2+(q-1)R\}\{R^2+(q-1)R\}+\left(-4R+4-2(q-1)\right)\{R^2+(q-1)R\}
\end{gather*}
and $\left(r^{-2}\{-R^2-(q-1)R\}\right)(r^{-2}\triangle_{S^q})$ equals
\[
r^{-4}\triangle_{S^q}\left(\{-R^2-(q-1)R\}+4R-4+2(q-1)\right) ,
\]
on $S^p\times S^q$,
\begin{gather*}
 \dA^2 =\triangle_{S^q}^2+\triangle_{S^p}^2-2\triangle_{S^q}\triangle_{S^p} \\
\phantom{\dA^2 =}{} +2\triangle_{S^q}\left(\{-R^2-(q-1)R\}+\{S^2+(p-1)S\}+2R-2+(q-1)\right) \\
\phantom{\dA^2 =}{} -2\triangle_{S^p}\left(\{-R^2-(q-1)R\}+\{S^2+(p-1)S\}-2S+2-(p-1)\right) \\
\phantom{\dA^2 =}{} +\{-R^2-(q-1)R\}^2+\{S^2+(p-1)S\}^2+2\{-R^2-(q-1)R\}\{S^2+(p-1)S\} \\
\phantom{\dA^2 =}{} +(4R-4+2(q-1))\{-R^2-(q-1)R\}+(-4S+4-2(p-1))\{S^2+(p-1)S\} .
\end{gather*}
Let
\begin{gather*}
A:=\{-R^2-(q-1)R\}, \qquad B:=\{S^2+(p-1)S\},\\
C:=2R-2+(q-1), \qquad D:=-2S+2-(p-1) .
\end{gather*}
Then,
\begin{gather*}
\dA^2 \mid_{S^p\times S^q}=\triangle_{S^q}^2+\triangle_{S^p}^2-2\triangle_{S^q}\triangle_{S^p}
      +2\triangle_{S^q}(A+B+C)-2\triangle_{S^p}(A+B+D) \\
\phantom{\dA^2 \mid_{S^p\times S^q}=}{} +(A+B+C)(A+B+D)+(A-B)(C-D)-CD .
\end{gather*}

Note that, on ${\cal E}[\omega]$,
\begin{gather*}
A+B+C=\left(-\omega- {\frac{n}{2}}+2\right)Y-\left( {\frac{q-p}{2}}-1\right)\omega +(q-3),\\
A+B+D=\left(-\omega- {\frac{n}{2}}+2\right)Y-\left( {\frac{q-p}{2}}+1\right)\omega -(p-3) .
\end{gather*}

Note also that
\begin{gather*}
(A-B)(C-D)-CD=Y^2\left(-\omega- {\frac{n}{2}}+2\right)
+Y\left\{(q-p)\left(-\omega- {\frac{n}{2}}+2\right)\right\}\\
\phantom{(A-B)(C-D)-CD=}{} +\left\{\omega^2+n\omega-2\omega+(\omega+q-3)(\omega+p-3)\right\} .
\end{gather*}
Since $\omega=2-\frac{n}{2}$,
\begin{gather*}
\left( {\frac{q-p}{2}}-1\right)\omega -(q-3)=1- {\frac{n-2}{4(n-1)}}\mbox{Scal} ,\\
\left( {\frac{q-p}{2}}+1\right)\omega +(p-3)=-1- {\frac{n-2}{4(n-1)}}\mbox{Scal},\qquad \mbox{and}\\
\omega^2+n\omega-2\omega+(\omega+q-3)(\omega+p-3)=n- {\frac{n^2}{2}}+pq+1,
\end{gather*}
we have
\begin{gather*}
\dA^2 \mid_{S^p\times S^q}=\triangle_{S^q}^2+\triangle_{S^p}^2-2\triangle_{S^q}\triangle_{S^p}
       -2\left\{1- {\frac{n-2}{4(n-1)}}\mbox{Scal}\right\}\triangle_{S^q}
\\
\phantom{\dA^2 \mid_{S^p\times S^q}=}{} +2\left\{-1- {\frac{n-2}{4(n-1)}}\mbox{Scal}\right\}\triangle_{S^p}
+\left( {\frac{n-2}{4(n-1)}}\mbox{Scal}\right)^2+n- {\frac{n^2}{2}}+pq .
\end{gather*}
We claim that this is the Paneitz operator \cite{Branson:95}
\[
P=\triangle^2+\delta T d+ {\frac{n-4}{2}}Q ,
\]
where
\begin{gather*}
J=\mbox{Scal}/(2(n-1)) , \qquad
V=(\rho-Jg)/(n-2) , \\
T=(n-2)J-4V\cdot, \qquad
Q= {\frac{n}{2}}J^2-2|V|^2+\triangle J .
\end{gather*}
Since ${\rm  Scal}=(n-1)(q-p)$ and $J= {\frac{q-p}{2}}$,
\begin{gather*}
\delta J d= {\frac{q-p}{2}}(\triangle_{S^q}-\triangle_{S^p}) ,
\\
\delta V d= {\frac{1}{n-2}}\left\{(p-1)\triangle_{S^p}+(q-1)\triangle_{S^q}\right\}
- {\frac{q-p}{2(n-2)}}(\triangle_{S^q}-\triangle_{S^p}) .
\end{gather*}

Thus
\begin{gather*}
\delta T d= {(n-2)\frac{q-p}{2}(\triangle_{S^q}-\triangle_{S^p})}
 +\frac{1}{n-2}\{-4(p-1)-2(q-p)\}\triangle_{S^p}\\
\phantom{\delta T d=}{} +
	  \frac{1}{n-2}\{-4(q-1)+2(q-p)\}\triangle_{S^q}\\
\phantom{\delta T d}{} = {2\left(\frac{n-2}{4(n-1)}\mbox{Scal}\right)(\triangle_{S^q}-\triangle_{S^p})
	  -2(\triangle_{S^q}+\triangle_{S^p})} .
\end{gather*}
On the other hand, since $|V|^2= {\frac{n}{4}}$,
\begin{gather*}
\frac{n-4}{2}Q= {\frac{n-4}{2}}\left( {\frac{n}{2}}
\left( {\frac{q-p}{2}}\right)^2- {\frac{n}{2}}\right)
              = {\frac{n(n-4)}{(n-2)^2}}\left( {\frac{n-2}{4(n-1)}}
              \mbox{Scal}\right)^2- {\frac{n(n-4)}{4}}\\
\phantom{\frac{n-4}{2}Q}{} =\left( {\frac{n-2}{4(n-1)}}\mbox{Scal}\right)^2-
       {\frac{4}{(n-2)^2}}\left( {\frac{n-2}{4(n-1)}}
      \mbox{Scal}\right)^2- {\frac{n(n-4)}{4}}\\
\phantom{\frac{n-4}{2}Q}{}=\left( {\frac{n-2}{4(n-1)}}\mbox{Scal}\right)^2+n
- {\frac{n^2}{2}}+pq
\end{gather*}
and the claim follows.

\section{Higher order operators}

Let
\[
C:=\sqrt{\Delta_{S^q}+\left( {\frac{q-1}2}\right)^2},\qquad
B:=\sqrt{\Delta_{S^p}+\left( {\frac{p-1}2}\right)^2},
\]
so that $C$ and $B$ are nonnegative operators with
\[
\Delta_{S^q}=C^2-\left(\frac{q-1}2\right)^2,\qquad
\Delta_{S^p}=B^2-\left(\frac{p-1}2\right)^2.
\]
The eigenvalue list for $\triangle_{S^q}$ \cite{IT:78, Branson:92} is
\[
j(q-1+j),\qquad j=0,1,2,\dots,
\]
so the eigenvalue list for $C$ is
\begin{gather}\label{eiglistT}
j+ {\frac{q-1}2},\qquad j=0,1,2,\dots.
\end{gather}
Similarly, the eigenvalue list for $B$ is
\begin{gather}\label{eiglistH}
k+ {\frac{p-1}2},\qquad k=0,1,2,\dots.
\end{gather}
Applying $\dA^m$, we get (with $k=m-\ell$)
\begin{gather*}
\dA^mf= \sum_{\ell=0}^m(-1)^k\binom{m}{\ell}s^{-2\ell}r^{-2k}
\left(C+\frac{q-1}2+v\right)
\cdots\left(C+\frac{q-1}2-2(\ell-1)+v\right)\bullet \\
\phantom{\dA^mf=}{} \left(C-\frac{q-1}2-v\right)
\cdots\left(C-\frac{q-1}2+2(\ell-1)-v\right) \bullet \\
\phantom{\dA^mf=}{} \left(B+\frac{p-1}2+u\right)
 \cdots\left(B+\frac{p-1}2-2(k-1)+u\right) \bullet \\
\phantom{\dA^mf=}{} \left(B-\frac{p-1}2-u\right)
\cdots\left(B-\frac{p-1}2+2(k-1)-u\right)f,
\end{gather*}
where in each $\cdot$, we move in increments or decrements of $2$.
These increments and decrements are determined by the homogeneity
drops implemented by the $s^{-2}$ and $r^{-2}$ factors in \eqref{dAlem}.
To restrict to $S^p\times S^q$, we just set $s=r=1$.

As a result, with
\[
Q:=\frac{q-1}2+v,\qquad P:=\frac{p-1}2+u,
\]
as long as we have the correct weight condition (equivalent to
(\ref{correctw}))
\[
P+Q=m-1,
\]
then the operator
\begin{gather*}
A_{2m}(C,B,Q):=\sum_{\ell=0}^m(-1)^k\binom{m}{\ell}\left(C+Q\right)
\cdots\left(C+Q-2(\ell-1)\right) \bullet \\
\phantom{A_{2m}(C,B,Q):=}{} \left(C-Q\right)
\cdots\left(C-Q+2(\ell-1)\right) \bullet\\
\phantom{A_{2m}(C,B,Q):=}{} \left(B+P\right)
\cdots\left(B+P-2(k-1)\right) \bullet \\
\phantom{A_{2m}(C,B,Q):=}{}  \left(B-P\right)
\cdots\left(B-P+2(k-1)\right), \qquad
k+\ell=m, \qquad
P+Q=m-1,
\end{gather*}
intrinsically def\/ined on $S^p\times S^q$,
intertwines $u_{m-n/2}$ and $u_{-m-n/2}$.

Note that the dependence of $A_{2m}(C,B,Q)$ is only on $(m,C,B,Q)$, since
$(u,P)$ is determined by $(m,C,B,Q)$.  The notation suggests substituting
numerical values for $C$ and $B$, a procedure justif\/ied by the eigenvalue lists
(\ref{eiglistT}), (\ref{eiglistH}).  These numerical values are nonnegative
real numbers, and depending on the parities of $q$ and $p$, they
are either integral or properly half-integral.

We {\it claim} that
\begin{proposition}
\begin{gather}
A_{2m}(C,B,Q)= (C+B+m-1)\cdots(C+B-m+1)\nonumber\\
\phantom{A_{2m}(C,B,Q)=}{}\times (C-B+m-1)\cdots(C-B-m+1)
:=G_{2m}(C,B),\label{theclaim}
\end{gather}
where the decrements are by $2$ units each time.
\end{proposition}
In particular, we are claiming that the left-hand side of (\ref{theclaim}) is independent of $Q$.

The operator $G_{2m}(C,B)$ is in fact a dif\/ferential operator since
\begin{gather*}
\left\{\!
\begin{array}{l}
\displaystyle (C+B)(C-B)  {\prod_{l=1}^{(m-1)/2}}\left[C+(B+2l)\right]\left[C-(B+2l)\right]\vspace{1mm}\\
\qquad{} \times \left[C+(B-2l)\right]\left[C-(B-2l)\right], \qquad m\mbox{ odd},\vspace{1mm}\\
\displaystyle {\prod_{l=1}^{m/2}}\left[C+\left(B+(2l-1)\right)\right]\left[C-\left(B+(2l-1)\right)\right]\vspace{1mm}\\
\qquad{} \times \left[C+\left(B-(2l-1)\right)\right]\left[C-\left(B-(2l-1)\right)\right], \qquad m\mbox{ even},
\end{array}\right.
\\
\qquad{}=\left\{\!
\begin{array}{l}
\displaystyle (C^2-B^2)
 {\prod_{l=1}^{(m-1)/2}}
\left[C^4-2(B^2+(2l)^2)C^2+(B^2-(2l)^2)^2\right],
\qquad m\mbox{ odd},\vspace{1mm}\\
\displaystyle {\prod_{l=1}^{m/2}}
\left[C^4-2(B^2+(2l-1)^2)C^2+(B^2-(2l-1)^2)^2\right],
\qquad m\mbox{ even}.
\end{array}\right.
\end{gather*}

Recently, Gover \cite{Gover:06} showed that on conformally Einstein manifolds, the operators are of the form
\[
\square_m=\prod_{l=1}^{m}(\Delta-c_l\mbox{Sc}) ,
\]
where $c_l=(n+2l-2)(n-2l)/(4n(n-1))$, Sc is the scalar curvature and $\Delta=\nabla^a\nabla_a$.

To get the formula (\ref{theclaim}) in case of sphere $S^n$, we set
\[
C:=\sqrt{\Delta+\left( {\frac{n-1}2}\right)^2}, \qquad B:=\frac{1}{2}.
\]
And the formula simplif\/ies to
\begin{gather*}
G_{2m}(C,1/2)= {\prod_{l=1}^{m}}\left(C- {\frac{2l-1}{2}}\right)\left(C+ {\frac{2l-1}{2}}\right)\\
\phantom{G_{2m}(C,1/2)}{}  = {\prod_{l=1}^{m}}\left(\Delta+\underbrace{ {\frac{(n+2l-2)(n-2l)}{4n(n-1)}}}_{c_l}
\underbrace{n(n-1)}_{\mbox{Sc}}\right) .
\end{gather*}
The ``+'' sign in the above is due to our convention $\Delta=-\nabla^a\nabla_a$ so the two formulas $\square_m$ and $G_{2m}(C,1/2)$ agree.

We will now prove the equality in (\ref{theclaim}).

Because of the eigenvalue lists (\ref{eiglistT}), (\ref{eiglistH}), to prove this in the
case in which $q$ and $p$ are odd, it is suf\/f\/icient to prove the identity
(\ref{theclaim}) with $C$ and $B$ replaced by nonnegative integers.  This will hold, in
turn, if it holds for $q=p=1$, the explicit mention of the dimensions having disappeared
in (\ref{theclaim}).

To prove (\ref{theclaim}) for $q=p=1$, note that each expression is polynomial in
$(C,B,v)$ for f\/ixed $m$, and that the highest degree terms in $(C,B)$
add up to $(C^2-B^2)^m$ for each expression.  Thus it will be enough to prove that
the right-hand side of (\ref{theclaim}) is the unique (up to constant multiples)
intertwinor $u_{m-1}\to u_{-m-1}$
in the case $q=p=1$.

By $K=$SO$(2)\times$SO$(2)$ invariance, an intertwinor $A$ must take an eigenvalue on
each
\[
\varphi_{j,f}:=e^{\sqrt{-1} ft}e^{\sqrt{-1} j\rho}
\]
for $\rho$ and $t$ the usual angular parameters on the positive-metric $S^1$ and the negative-metric $S^1$
respectively, and $f$ and $j$ integers.

The prototypical conformal vector f\/ield \cite{Branson:87} is
\[
T=\cos(\rho)\sin(t)\partial_t+\cos(t)\sin(\rho)\partial_\rho,
\]
with conformal factor
\[
\omega=\cos(\rho)\cos(t).
\]
The representation $U_{-r}$ of the Lie algebra $\mathfrak{so}(2,2)$ has
\begin{gather}
U_{-r}(T)\varphi_{j,f}=\frac14\big\{(f+j+r)\varphi_{j+1,f+1}+(f-j+r)\varphi_{j-1,f+1} \nonumber\\
\phantom{U_{-r}(T)\varphi_{j,f}=}{} +(-f+j+r)\varphi_{j+1,f-1}+(-f-j+r)\varphi_{j-1,f-1}\big\}.\label{whatUrdoes}
\end{gather}

Consider the operator
\[
P(\varepsilon,\delta):\varphi_{j,f}\mapsto\varphi_{j+\varepsilon,f+\delta}
\]
for $\varepsilon,\delta\in\{\pm 1\}$,
and the operators
\[
J:\varphi_{j,f}\mapsto j\varphi_{j,f},\qquad F:\varphi_{j,f}\mapsto f\varphi_{j,f}.
\]

Note that another expression for $G_{2m}(C,B)$ is as $G_{2m}(J,F)$,
since
\[
G_{2m}(J,F)=G_{2m}(J,-F)=G_{2m}(-J,F) .
\]
We have
\[
JP(\varepsilon,\delta)=P(\varepsilon,\delta)(J+\varepsilon),\qquad FP(\varepsilon,\delta)=P(\varepsilon,\delta)(F+\delta).
\]
By (\ref{whatUrdoes}),
\begin{gather*}
U_{\pm m-1}(T)=\frac14\big\{P(1,1)(J+F+1\mp m)+P(-1,1)(-J+F+1\mp m) \\
\phantom{U_{\pm m-1}(T)=}{} +P(1,-1)(J-F+1\mp m)+P(-1,-1)(-J-F+1\mp m)\big\}.
\end{gather*}
With this we may compute that
\begin{gather*}
4G_{2m}(J,F)U_{m-1}(T)=
G_{2m}(J,F)\big\{P(1,1)(J+F+1-m)+P(-1,1)(-J+F+1-m) \\
\phantom{4G_{2m}(J,F)U_{m-1}(T)=}{}+P(1,-1)(J-F+1-m)+P(-1,-1)(-J-F+1-m)\big\} \\
\phantom{4G_{2m}(J,F)U_{m-1}(T)}{} =P(1,1)\big\{(J+F+m+1)\cdots(J+F-m+3)\bullet \\
\phantom{4G_{2m}(J,F)U_{m-1}(T)=}{} (J-F+m-1)\cdots(J-F-m+1)\big\}
(J+F-m+1) \\
\phantom{4G_{2m}(J,F)U_{m-1}(T)=}{}+P(-1,1)\big\{(J+F+m-1)\cdots(J+F-m+1)\bullet \\
\phantom{4G_{2m}(J,F)U_{m-1}(T)=}{}(J-F+m-3)\cdots(J-F-m-1)\big\}
(-J+F-m+1) \\
\phantom{4G_{2m}(J,F)U_{m-1}(T)=}{}+P(1,-1)\big\{(J+F+m-1)\cdots(J+F-m+1)\bullet \\
\phantom{4G_{2m}(J,F)U_{m-1}(T)=}{}(J-F+m+1)\cdots(J-F-m+3)\big\}
(J-F-m+1) \\
\phantom{4G_{2m}(J,F)U_{m-1}(T)=}{}+P(-1,-1)\big\{(J+F+m-3)\cdots(J+F-m-1)\bullet \\
\phantom{4G_{2m}(J,F)U_{m-1}(T)=}{} (J-F+m-1)\cdots(J-F-m+1)\big\}
(-J-F-m+1),
\end{gather*}
whereas
\begin{gather*}
4U_{-m-1}(T)G_{2m}(J,F)=  \big\{P(1,1)(J+F+m+1)+P(-1,1)(-J+F+m+1) \\
\qquad{} +P(1,-1)(J-F+m+1) +P(-1,-1)(-J-F+m+1)\big\}\bullet \\
\qquad{} (J+F+m-1)\cdots(J+F-m+1)(J-F+m-1)\cdots(J-F-m+1).
\end{gather*}
The right-hand sides of the two preceding displays agree, so we have an intertwining
operator.

As a corollary, the claim (\ref{theclaim}) follows, so that $G_{2m}(C,B)$ is an
intertwinor whenever $qp$ is odd.  In fact, by polynomial continuation
from positive integral values,
the identity (\ref{theclaim}) holds whenever any complex values are substituted
for $C$ and $B$.  In particular, we can substitute proper half-integers,
and thus remove the condition that $qp$ be odd.

It would be good to have a proof which avoids a dimensional continuation argument. We present in the following appendix a proof which uses only an elementary combinatorial argument.

\appendix

\section{Appendix}
Here we use induction on the order of the operator. We will do:
\begin{itemize}\itemsep=0pt
\item Express $A_{2(m+1)}(C,B,Q)$ in terms of $(Q-1)$ in all terms containing $B$.
\item Compute and see
\begin{gather*}
A_{2m}(C-1,B,Q-1) \{A_2(C+m,B,Q+m)+A_2(C+m,B,Q-m)\}\\
\qquad{}=A_{2m}(C-1,B,Q-1)(C+m-B)(C+m+B) + A_{2(m+1)}(C,B,Q).
\end{gather*}
\item Since the above simply says
\begin{gather*}
2  A_{2m}(C-1,B,Q-1)  (C+m-B)(C+m+B)\\
\qquad{}=A_{2m}(C-1,B,Q-1)  (C+m-B)(C+m+B)
+A_{2(m+1)}(C,B,Q) ,
\end{gather*}
conclude
\[
A_{2(m+1)}(C,B,Q)=A_{2m}(C-1,B,Q-1) (C+m-B)(C+m+B).
\]
\end{itemize}

To go on, note f\/irst that $A_{2(m+1)}(C,B,Q)$ in terms of $(Q-1)$ in all terms containing $B$ is
\begin{gather}
 \sum_{\ell=0}^{m+1}(-1)^{m+1-\ell}\binom{m+1}{\ell}\nonumber\\
\qquad {}\bullet \left(B-(Q-1)+(m-1)\right)
\cdots\left(B-(Q-1)-(m-1)+2(\ell-1)\right)\nonumber\\
\qquad {}\bullet \left(B+(Q-1)-(m-1)\right)
\cdots\left(B+(Q-1)+(m-1)-2(\ell-1)\right)\nonumber\\
\qquad {} \bullet\left(C-Q\right)
\cdots\left(C-Q+2(\ell-1)\right)\nonumber\\
\qquad {} \bullet \left(C+Q\right)
\cdots\left(C+Q-2(\ell-1)\right) .\label{m+1}
\end{gather}

$A_{2m}(C-1,B,Q-1)$ can be written
\begin{gather}
 \sum_{\ell=0}^m(-1)^{m-\ell}\binom{m}{\ell}\nonumber\\
\qquad {}\bullet \left(B-(Q-1)+(m-1)\right)
\cdots\left(B-(Q-1)-(m-1)+2\ell\right)\nonumber\\
\qquad {}\bullet \left(B+(Q-1)-(m-1)\right)
\cdots\left(B+(Q-1)+(m-1)-2\ell\right)\nonumber\\
\qquad {}\bullet\left(C-Q\right)
\cdots\left(C-Q+2(\ell-1)\right)\nonumber\\
\qquad {}\bullet \left(C+Q-2\right)
\cdots\left(C+Q-2\ell\right).\label{m1}
\end{gather}
Def\/ine $R_B$ and $R_C$ to express $A_{2m}(C-1,B,Q-1)$
as
\begin{gather*}
(-1)^{m}\binom{m}{0}
 \bullet \left(B-(Q-1)+(m-1)\right)
\cdots\left(B-(Q-1)-(m-1)\right)\\
\qquad{}\bullet \left(B+(Q-1)-(m-1)\right)
\cdots\left(B+(Q-1)+(m-1)\right)+R_B
\end{gather*}
or
\begin{gather*}
R_C +(-1)^{0}\binom{m}{m}
\bullet\left(C-Q\right)
\cdots\left(C-Q+2(m-1)\right)\\
\qquad{} \bullet \left(C+Q-2\right)
\cdots\left(C+Q-2m\right) .
\end{gather*}
We can write
\begin{gather*}
A_{2}(C+m,B,Q+m)=(C+m+B)(C+m-B)\\
\phantom{A_{2}(C+m,B,Q+m)}{}  = -(B-(Q+m))(B+(Q+m))+ (C-Q)(C+Q+2m)
\end{gather*}
and
\begin{gather*}
A_{2}(C+m,B,Q-m)=(C+m+B)(C+m-B)\\
\phantom{A_{2}(C+m,B,Q-m)}{}  =-(B-(Q-m))(B+(Q-m))
              +(C-Q+2m)(C+Q) .
\end{gather*}
The f\/irst product $A_{2m}(C-1,B,Q-1) A_{2}(C+m,B,Q+m)$
becomes
\begin{gather*}
(-1)^{m}\binom{m}{0}\bullet \left(B-(Q-1)+(m-1)\right)
\cdots\left(B-(Q-1)-(m-1)\right)\\
\qquad {} \bullet \left(B+(Q-1)-(m-1)\right)
\cdots\left(B+(Q-1)+(m-1)\right)\\
\qquad{} \bullet \left\{-(B-(Q+m))(B+(Q+m))+(C-Q)(C+Q+2m)
\right\}\\
\qquad{} +
R_B\bullet (C+m+B)(C+m-B) ,
\end{gather*}
which can be rewritten as
\begin{gather*}
(-1)^{m+1}\binom{m+1}{0}
\bullet \left(B-(Q-1)+(m-1)\right)
\cdots\left(B-(Q-1)-(m-1)\right)\left(B-(Q+m)\right)\\
\qquad {} \bullet \left(B+(Q-1)-(m-1)\right)
\cdots\left(B+(Q-1)+(m-1)\right)\left(B+(Q+m)\right)\\
\qquad{}+
(-1)^{m}\binom{m}{0} \bullet \left(B-(Q-1)+(m-1)\right)
\cdots\left(B-(Q-1)-(m-1)\right)\\
\qquad{} \bullet \left(B+(Q-1)-(m-1)\right)
\cdots\left(B+(Q-1)+(m-1)\right)\\
\qquad{} \bullet (C-Q)(C+Q+2m)
+
R_B\bullet (C+m+B)(C+m-B) .
\end{gather*}
The second product $A_{2m}(C-1,B,Q-1)  A_{2}(C+m,B,Q-m)$
is
\begin{gather*}
R_C\bullet (C+m+B)(C+m-B)
+(-1)^{0}\binom{m}{m}\bullet\left(C-Q\right)
\cdots\left(C-Q+2(m-1)\right)\\
\qquad{} \bullet \left(C+Q-2\right)
\cdots\left(C+Q-2m\right)\\
\qquad{} \bullet \left\{-(B-(Q-m))(B+(Q-m))+(C-Q+2m)(C+Q)
\right\}\\
\qquad{} = R_C\bullet (C+m+B)(C+m-B)
+
(-1)^{1}\binom{m}{m}\bullet\left(C-Q\right)
\cdots\left(C-Q+2(m-1)\right)\\
\qquad{} \bullet \left(C+Q-2\right)
\cdots\left(C+Q-2m\right)\bullet (B-(Q-m))(B+(Q-m))\\
\qquad{} +
(-1)^{0}\binom{m+1}{m+1}
\bullet\left(C-Q\right)
\cdots\left(C-Q+2m\right)
\bullet \left(C+Q\right)
\cdots\left(C+Q-2m\right)  .
\end{gather*}
So by adding up the above two products, we get
\begin{gather*}
(-1)^{m+1}\binom{m+1}{0}
\bullet \left(B-(Q-1)+(m-1)\right)
\cdots\left(B-(Q-1)-(m-1)\right)\left(B-(Q+m)\right)\\
\qquad {} \bullet \left(B+(Q-1)-(m-1)\right)
\cdots\left(B+(Q-1)+(m-1)\right)\left(B+(Q+m)\right)\\
\qquad{} +
(-1)^{0}\binom{m+1}{m+1}
\bullet\left(C-Q\right)
\cdots\left(C-Q+2m\right)
\bullet \left(C+Q\right)
\cdots\left(C+Q-2m\right)\\
\qquad{} +
(-1)^{m}\binom{m}{0} \bullet \left(B-(Q-1)+(m-1)\right)
\cdots\left(B-(Q-1)-(m-1)\right)\\
\qquad{} \bullet \left(B+(Q-1)-(m-1)\right)
\cdots\left(B+(Q-1)+(m-1)\right)\bullet (C-Q)(C+Q+2m)\\
\qquad{} +
(-1)^{1}\binom{m}{m}
\bullet\left(C-Q\right)
\cdots\left(C-Q+2(m-1)\right)
\bullet \left(C+Q-2\right)
\cdots\left(C+Q-2m\right)\\
\qquad{} \bullet (B-(Q-m))(B+(Q-m))\\
\qquad{} +
(R_B+R_C)\bullet (C+m+B)(C+m-B) .
\end{gather*}

Note that $R_B$ is missing the f\/irst term  and $R_C$ is
missing the last term of (\ref{m1}). So we have
\begin{gather*}
(R_B+R_C)\bullet (C+m+B)(C+m-B)\\
\qquad {}=
A_{2m}(C-1,B,Q-1)  (C+m+B)(C+m-B)\\
\qquad{} +  \sum_{\ell=1}^{m-1}(-1)^{m-\ell}\binom{m}{\ell}
\bullet \left(B-(Q-1)+(m-1)\right)
\cdots\left(B-(Q-1)-(m-1)+2\ell\right)\\
\qquad {} \bullet \left(B+(Q-1)-(m-1)\right)
\cdots\left(B+(Q-1)+(m-1)-2\ell\right)\\
\qquad{} \bullet\left(C-Q\right)
\cdots\left(C-Q+2(\ell-1)\right)
\bullet \left(C+Q-2\right)
\cdots\left(C+Q-2\ell\right)\\
\qquad{} \bullet (C+m+B)(C+m-B) .
\end{gather*}
Therefore, $A_{2m}(C-1,B,Q-1)  (C+m+B)(C+m-B)$ equals (see (\ref{m+1}))
\begin{gather}
1\mbox{st term of }A_{2(m+1)}(C,B)+(m+2)\mbox{nd term of }A_{2(m+1)}(C,B)
\nonumber\\[3mm]
\quad {} + (-1)^{m}\binom{m}{0}
\bullet \left(B-(Q-1)+(m-1)\right)
\cdots\left(B-(Q-1)-(m-1)\right)\nonumber\\
\qquad {} \bullet \left(B+(Q-1)-(m-1)\right)
\cdots\left(B+(Q-1)+(m-1)\right)\nonumber\\
\qquad{} \bullet (C-Q)(C+Q+2m)\label{3rd}\\[3mm]
\quad{} +
(-1)^{1}\binom{m}{m}\bullet\left(C-Q\right)
\cdots\left(C-Q+2(m-1)\right)\nonumber\\
\qquad{}\bullet \left(C+Q-2\right)
\cdots\left(C+Q-2m\right)
\bullet (B-(Q-m))(B+(Q-m))\label{4th}\\[3mm]
\quad{}+
 \sum_{\ell=1}^{m-1}(-1)^{m-\ell}\binom{m}{\ell}
\bullet \left(B-(Q-1)+(m-1)\right)
\cdots\left(B-(Q-1)-(m-1)+2\ell\right)\nonumber\\
\qquad{} \bullet \left(B+(Q-1)-(m-1)\right)
\cdots\left(B+(Q-1)+(m-1)-2\ell\right)\nonumber\\
\qquad {} \bullet\left(C-Q\right)
\cdots\left(C-Q+2(\ell-1)\right)
\bullet \left(C+Q-2\right)
\cdots\left(C+Q-2\ell\right)\nonumber\\
\qquad{} \bullet (C+m+B)(C+m-B) .\label{5th}
\end{gather}
So we want to show (\ref{3rd})\;+\;(\ref{4th})\;+\;(\ref{5th}) is exactly
the other term in $A_{2(m+1)}(C,B)$.

Since, for any $Q_{\ell}$,
\[
(C+m+B)(C+m-B)=\{-(B-Q_{\ell})(B+Q_{\ell})
+(C+m-Q_{\ell})(C+m+Q_{\ell})\} ,
\]
(\ref{5th}) becomes
\begin{gather*}
 \sum_{\ell=1}^{m-1}(-1)^{m+1-\ell}\binom{m}{\ell}
\bullet \left(B-(Q-1)+(m-1)\right)
\cdots\left(B-(Q-1)-(m-1)+2\ell\right)\\
\qquad{} \bullet \left(B+(Q-1)-(m-1)\right)
\cdots\left(B+(Q-1)+(m-1)-2\ell\right)\\
\qquad{} \bullet\left(C-Q\right)
\cdots\left(C-Q+2(\ell-1)\right)
\bullet \left(C+Q-2\right)
\cdots\left(C+Q-2\ell\right)
\bullet (B-Q_{\ell})(B+Q_{\ell})\\
\qquad{} +
 \sum_{\ell=1}^{m-1}(-1)^{m-\ell}\binom{m}{\ell}
\bullet \left(B-(Q-1)+(m-1)\right)
\cdots\left(B-(Q-1)-(m-1)+2\ell\right)\\
\qquad{} \bullet \left(B+(Q-1)-(m-1)\right)
\cdots\left(B+(Q-1)+(m-1)-2\ell\right)\\
\qquad{} \bullet\left(C-Q\right)
\cdots\left(C-Q+2(\ell-1)\right)
\bullet \left(C+Q-2\right)
\cdots\left(C+Q-2\ell\right)\\
\qquad{} \bullet (C+m-Q_{\ell})(C+m+Q_{\ell})\\
\qquad{} =:
 \sum_{\ell=1}^{m-1}(-1)^{m+1-\ell}\binom{m}{\ell}
E_{\ell}
+
 \sum_{\ell=1}^{m-1}(-1)^{m-\ell}\binom{m}{\ell}
F_{\ell}
\\
\qquad{}=
 \sum_{\ell=2}^{m-2}(-1)^{m+1-\ell}
\left(\binom{m}{\ell}E_{\ell}
+\binom{m}{\ell-1}F_{\ell-1}\right)\\
\qquad{} +(-1)^{m}\binom{m}{1}E_1+(-1)^{1}\binom{m}{m-1}F_{m-1} .
\end{gather*}
But $\left(\binom{m}{\ell}E_{\ell}
+\binom{m}{\ell-1}F_{\ell-1}\right)$ becomes
\begin{gather*}
\binom{m}{\ell}
 \left(B-(Q-1)+(m-1)\right)
\cdots\left(B-(Q-1)-(m-1)+2(\ell-1)\right)\\
\qquad{} \bullet \left(B+(Q-1)-(m-1)\right)
\cdots\left(B+(Q-1)+(m-1)-2(\ell-1)\right)\\
\qquad{} \bullet\left(C-Q\right)
\cdots\left(C-Q+2(\ell-1)\right)
\bullet \left(C+Q-2\right)
\cdots\left(C+Q-2\ell\right)\\
\qquad{} +
\binom{m}{\ell-1}
 \left(B-(Q-1)+(m-1)\right)
\cdots\left(B-(Q-1)-(m-1)+2(\ell-1)\right)\\
\qquad{} \bullet \left(B+(Q-1)-(m-1)\right)
\cdots\left(B+(Q-1)+(m-1)-2(\ell-1)\right)\\
\qquad{} \bullet\left(C-Q\right)
\cdots\left(C-Q+2((\ell-1)-1)\right)
\bullet \left(C+Q-2\right)
\cdots\left(C+Q-2(\ell-1)\right)\\
\qquad{} \bullet (C+m-Q_{\ell-1})(C+m+Q_{\ell-1})
\end{gather*}
which is, upon choosing $Q_{\ell}$ to be $Q-2\ell+m$,
\begin{gather*}
 \left(B-(Q-1)+(m-1)\right)
\cdots\left(B-(Q-1)-(m-1)+2(\ell-1)\right)\\
\qquad{} \bullet \left(B+(Q-1)-(m-1)\right)
\cdots\left(B+(Q-1)+(m-1)-2(\ell-1)\right)\\
\qquad{} \bullet\left(C-Q\right)
\cdots\left(C-Q+2(\ell-1)\right)
\bullet \left(C+Q-2\right)
\cdots\left(C+Q-2(\ell-1)\right)
\end{gather*}
times
\[
\binom{m}{\ell}(C+Q-2\ell)+\binom{m}{\ell-1}(C+Q-2(\ell-1-m)) ,
\]
since
\begin{gather*}
B-Q_{\ell}=B-(Q-1)-(m-1)+2(\ell-1) ,\\
B+Q_{\ell}=B+(Q-1)+(m-1)-2(\ell-1) ,\\
C+m-Q_{\ell}=C-Q+2(\ell-1)\qquad
\mbox{and}\\
C+m+Q_{\ell}=C+Q+2(\ell-1-m) .
\end{gather*}

Note also that
\[
\binom{m}{\ell}(C+Q-2\ell)+\binom{m}{\ell-1}(C+Q-2(\ell-1-m))
=\binom{m+1}{l}(C+Q)  .
\]
Thus
\[
\left(\binom{m}{\ell}E_{\ell}
+\binom{m}{\ell-1}F_{\ell-1}\right)
=\mbox{ the }(\ell+1) \mbox{st term of }A_{2(m+1)}(C,B,Q) .
\]
Finally we note that (\ref{3rd}) and
$(-1)^{m}\binom{m}{1}E_1$ add up to
\begin{gather*}
(-1)^m
 \left(B-(Q-1)+(m-1)\right)
\cdots\left(B-(Q-1)-(m-1)\right)\\
\qquad{}\bullet \left(B+(Q-1)-(m-1)\right)
\cdots\left(B+(Q-1)+(m-1)\right)
\bullet (C-Q)
\end{gather*}
times
\[
\binom{m}{0}(C+Q+2m)+\binom{m}{1}(C+Q-2) ,
\]
since
\[
B-Q_1=B-(Q-1)-(m-1) \qquad \mbox{and} \qquad B+Q_1=B+(Q-1)+(m-1) .
\]
This is the 2nd term of $A_{2(m+1)}(C,B,Q)$, since
\[
\binom{m}{0}(C+Q+2m)+\binom{m}{1}(C+Q-2)
= \binom{m+1}{1}(C+Q)
\]
Similarly, (\ref{4th}) and
$(-1)^{1}\binom{m}{m-1}F_{m-1}$ add up to
\begin{gather*}
(-1)^{m}
\left(C-Q\right)
\cdots\left(C-Q+2(m-1)\right)
\bullet \left(C+Q-2\right)
\cdots\left(C+Q-2(m-1)\right)\\
\qquad{} \bullet (B-(Q-1)+(m-1))(B+(Q-1)-(m-1))
\end{gather*}
times
\[
\binom{m}{m}(C+Q-2m)+\binom{m}{m-1}(C+Q+2) ,
\]
since
\begin{gather*}
C+m-Q_{m-1}=C-Q+2(m-1) ,\\
C+m+Q_{m-1}=C+Q+2 ,\\
B-(Q-m)=B-(Q-1)+(m-1)\qquad \mbox{and}\\
B+(Q-m)=B+(Q-1)-(m-1) .
\end{gather*}
This is the $(m+1)$st term of $A_{2(m+1)}(C,B,Q)$, since
\[
\binom{m}{m}(C+Q-2m)+\binom{m}{m-1}(C+Q+2)
=\binom{m+1}{m}(C+Q) .
\]

\subsection*{Acknowledgements}
The second author gratefully acknowledges support by the Eduard Cech Center for Algebra and Geometry Grant Nr. LC505
and support by the Seoul National University Grant Nr.~KRF 2005-070-C00007.

\pdfbookmark[1]{References}{ref}
\LastPageEnding

\end{document}